\documentclass[12pt]{article}
\usepackage{amsmath,amssymb}
\def\CONE{\mathcal{S}}
\def\Laplace{\mathcal{L}}
\def\M{\mathcal{M}}
\def\C{\mathbb{C}}
\def\R{\mathbb{R}}
\newtheorem{theorem}{\hspace*{\parindent}Theorem}

\newtheorem{corollary}{\hspace*{\parindent}Corollary}
\newcounter{remark}
\newcommand{\rem}{\par\refstepcounter{remark}\noindent\textbf{Remark \arabic{remark}.} }
\title{Generalized Stieltjes transforms: basic aspects}
\author{D.\,Karp\footnote{Institute of Applied Mathematics, Vladivostok, Russia,
e-mail:\,\emph{dimkrp@gmail.com}}~~and~E.\,Prilepkina\footnote{Institute
of Applied Mathematics, Vladivostok, Russia,
e-mail:\,\emph{pril-elena@yandex.ru}}}
\date{}
\begin{document}
\maketitle

\begin{center}
\parbox{12cm}{
\small\textbf{Abstract.}  The paper surveys the basic properties
of generalized Stieltjes functions including some new ones.  We
introduce the notion of the exact Stieltjes order and give a
criterion of exactness, simple sufficient conditions and some
prototypical examples. The paper includes an appendix, where we
define the left sided Riemann-Liouville and the right sided
Kober-Erdelyi fractional integrals of measures supported on half
axis and give inversion formulas for them. }
\end{center}

\bigskip

Keywords: \emph{generalized Stieltjes transform, generalized
Stieltjes function, fractional integral, fractional derivative,
exact Stieltjes order}

\bigskip

MSC2010: 26A48, 26A33, 44A15

\bigskip

\paragraph{1. Introduction.} The generalized Stieltjes transform of a non-negative measure $\mu$
supported on $[0,\infty)$ is defined by
$$
\int\limits_{[0,\infty)}\frac{\mu(du)}{(u+z)^{\alpha}},
$$
where $\alpha>0$ and we always choose the branch of the power
function which is positive on the positive half-axis. The measure
is assumed to produce a convergent integral for each
$z\in\C\!\setminus\!(-\infty,0]$ thus generating a function
holomorphic in $\C\!\setminus\!(-\infty,-r]$, where
$r=\inf\{x:~x\in\mathbf{supp}(\mu)\}$. Functions representable by
the above integral plus a non-negative constant are known as
generalized Stieltjes functions \cite{Sokal,Sumner},
\cite[Section~8]{WidderP}, \cite[Chapter~VIII]{Widder}.

The case $\alpha=1$ has been thoroughly studied by many authors
beginning with the classical work of Stieltjes \cite{Stieltjes}
followed by Krein (see \cite{KN} and references therein), Widder
\cite{WidderP,Widder}, Hirsch \cite{Hirsch}, Berg \cite{Berg} and
many others.  Among the most important tools facilitating such
study are the complex inversion formula due to Stieltjes
\cite{Stieltjes,Widder},  the complex variable characterization
found by Krein (see Theorem~\ref{th:Krein} below) and the real
inversion formulas by Widder \cite{WidderP,Widder}. When the
measure $\mu$ has compact support, the Stieltjes functions are
also known as Markov functions studied by Chebyshev and Markov in
connection with continued fractions.  The deep connection of the
Stieltjes and Markov functions with continued fractions and
Pad\'{e} approximation is investigated in the monographs
\cite{BG,CPVWJ}.  See also the survey paper \cite{Gilewicz}.
Connection with Bernstein functions and various other similar
classes can be found in the carefully written recent  monograph
\cite{SSV}.

For general $\alpha>0$ much less is known. A complex inversion
formula in this case has been found by Sumner \cite{Sumner} and
later rediscovered by Schwarz \cite{Schwarz}. It was amended by
several other complex inversion formulas by Byrne and Love in
\cite{Love1}.  These authors also found several real inversion
formulas in \cite{Love2,Love3}.  A simple real variable
characterization has been discovered recently by Sokal
\cite{Sokal} generalizing the corresponding result of Widder. The
asymptotic expansion of generalized Stieltjes transforms of slowly
decaying functions has been studied by L\'{o}pez and Ferreira in
\cite{LF}, factorization as iterative Laplace transforms - by
Y\"{u}rekli \cite{Yur}, closely related classes on half plane have
been investigated by Jerbashian \cite{Jerb}.  An interesting
connection to entire functions has been discovered in a recent
work of Pedersen \cite{Pedersen}.

In this paper we collect a number of facts about generalized
Stieltjes functions. Most of them are scattered in the literature,
those we could not find are furnished with detailed proofs. Some
of them may be new.  In the last section we introduce the notion
of the exact Stieltjes order, which is the  ''natural'' exponent
defined for each generalized Stieltjes function.  We give a
criterion of exactness and its two practical corollaries.  We also
added an appendix, where we define the left sided
Riemann-Liouville and the right sided Kober-Erdelyi fractional
integrals of measures supported on $[0,\infty]$ (the one point
compactification of $[0,\infty)$ - see details below) and give
inversion formulas for them.  Our study of the generalized
Stieltjes transforms carried out in this paper has been largely
motivated by the applications to the theory of hypergeometric
functions presented in our forthcoming work \cite{KaPr}.

\bigskip

\paragraph{2. Definition and real variable properties.}  Define $S_{\alpha}$,
$\alpha>0$, to be the class of functions representable by the
integral
\begin{equation}\label{eq:mu-rep}
f(z)=\int\limits_{(0,\infty)}\frac{\mu_{\alpha}(du)}{(u+z)^{\alpha}}+\mu_{\infty}+\frac{\mu_{0}}{z^{\alpha}},
\end{equation}
where $0\leq\mu_0,\mu_{\infty}<\infty$ and $\mu_{\alpha}$ runs
over the set of non-negative measures supported on $[0,\infty)$
such that
\begin{equation}\label{eq:mu-condition}
\int\limits_{[0,\infty)}\frac{\mu_{\alpha}(dt)}{(1+t)^{\alpha}}<\infty.
\end{equation}
This condition guarantees the finiteness of the integral
(\ref{eq:mu-rep})  for all $z\in\C\!\setminus\!(-\infty,0]$. The
classical Stieltjes cone $\CONE$ corresponds to the case
$\alpha=1$, i.e. $\CONE=S_1$, see \cite[formula (1)]{Berg},
\cite[formula (2)]{Sokal} or \cite[formula (2.1)]{SSV}.

If we define $[0,\infty]$ to be the one point compactification  of
$[0,\infty)$ formula (\ref{eq:mu-rep}) may be rewritten as
\begin{equation}\label{eq:mutilde-rep}
f(z)=\int\limits_{[0,\infty]}\left(\frac{u+1}{u+z}\right)^{\alpha}\tilde{\mu}_{\alpha}(du),
\end{equation}
where
$$
\tilde{\mu}_{\alpha}=\frac{\mu_{\alpha}(du)}{(1+u)^{\alpha}} +
\mu_{0}\delta_0+\mu_{\infty}\delta_{\infty},
$$
is a  finite measure on the compact interval $[0,\infty]$.  Here
$\delta_{a}$ stands for the Dirac measure with mass $1$
concentrated in $a$. This explains the notation $\mu_{\infty}$ for
the non-negative constant in (\ref{eq:mu-rep}). The set of
measures supported on $[0,\infty]$ and satisfying
(\ref{eq:mu-condition}) will be denoted by $\M_{\alpha}$. The
majority of references on the classical and generalized Stieltjes
functions use formula (\ref{eq:mu-rep}) (or its particular case
$\alpha=1$) to define them. See, for instance,
\cite{Berg,Hirsch,Love1,Love2,Love3,SSV,Sokal,Widder,Yur}.
However, the literature on Pad\'{e} approximation frequently
defines the Stieltjes functions by $\alpha=1$ case of the
following formula
\begin{equation}\label{eq:rho-rep}
f(z)=\int\limits_{(0,\infty)}\frac{\rho_{\alpha}(dt)}{(1+tz)^{\alpha}}+\rho_0+\frac{\rho_{\infty}}{z^{\alpha}}.
\end{equation}
Define the map $N_{\alpha}$ on $\M_{\alpha}$ by
\begin{equation}\label{eq:Nalpha-defined}
[N_{\alpha}\mu](A):=\int\limits_{1/A}t^{-\alpha}\mu(dt)=
\int\limits_{A}u^{\alpha}\mu^{\!*}(du)~\text{for each Borel
set}~A\subset(0,\infty),
\end{equation}
where  $\mu^{\!*}$ is the image measure of $\mu$ under the map
$t\to{1/t}$ and by definition
$[N_{\alpha}\mu](\{\infty\})=\mu(\{0\})$,
$[N_{\alpha}\mu](\{0\})=\mu(\{\infty\})$. If (\ref{eq:mu-rep}) and
(\ref{eq:rho-rep}) represent the same function, the change of
variable $t=1/u$ in (\ref{eq:mu-rep}) shows that
$\rho_{\alpha}=N_{\alpha}\mu_{\alpha}$.
 By the same
change of variable applied to (\ref{eq:mu-condition}) we see that
\begin{equation}\label{eq:rho-condition}
\int\limits_{[0,\infty)}\frac{\rho_{\alpha}(dt)}{(1+t)^{\alpha}}<\infty,
\end{equation}
so that $N_{\alpha}$ maps $\M_{\alpha}$ into itself. Moreover,
$N_{\alpha}$ is easily seen to be an involution on $\M_{\alpha}$:
$N_{\alpha}N_{\alpha}\mu=\mu$ for each $\mu\in\M_{\alpha}$.

Definition (\ref{eq:rho-rep}) is a natural extension  of the
definition of Stieltjes functions used in \cite[formula(5.1)]{BG}
and \cite[formula (1)]{Gilewicz}. In some situations this
representation leads to simpler expressions for hypergeometric
functions.  We will work with both representations
(\ref{eq:mu-rep}) and (\ref{eq:rho-rep}).  If we define the finite
measure
$$
\tilde{\rho}_{\alpha}=\frac{\rho_{\alpha}(du)}{(1+u)^{\alpha}} +
\rho_{0}\delta_0+\rho_{\infty}\delta_{\infty}
$$
on the compact interval $[0,\infty]$, then $\tilde{\rho}_{\alpha}$
and $\tilde{\mu}_{\alpha}$ from (\ref{eq:mutilde-rep}) are related
by $\tilde{\rho}_{\alpha}(A)=\tilde{\mu}_{\alpha}(1/A)$ for each
Borel set $A\subset[0,\infty]$.

We will denote by $F_{\mu}:[0,\infty)\to[0,\infty)$ the
left-continuous distribution function of the measure
$\mu\in\M_{\alpha}$: $F_{\mu}(x)=\mu([0,x))$ normalized by
$F_{\mu}(0)=0$. The distribution function $F_{\mu}$ defines the
measure $\mu$ uniquely except for a possible atom at infinity
which must be specified separately. According to
\cite[section~1.8]{Bogachev} or \cite[Problem~7.9(iii)]{Schilling}
every measure $\mu\in\M_{\alpha}$ is generated by such
non-decreasing left-continuous function satisfying
$$
\int\limits_{[0,\infty)}\frac{dF_{\mu}(t)}{(1+t)^{\alpha}}<\infty,
$$
and a non-negative constant $\mu_{\infty}:=\mu(\{\infty\})$. Here
the integral is understood as Lebesgue-Stieltjes integral
\cite[Section 20.3]{RF} or \cite[2.12(vi)]{Bogachev}.

The Stieltjes cone $\CONE$ possesses a number of nice stability
properties which can be found in \cite{Berg}.  The majority of
these properties do not carry over to $S_{\alpha}$,
$\alpha\ne{1}$. On the other hand, here we have some new effects
related to transition from $S_{\alpha}$ to $S_{\beta}$,
$\beta\ne{\alpha}$ and certain stability properties of the class
$S_{\infty}:=\cup_{\alpha>0}S_{\alpha}$. Below we list the basic
facts about these classes.

\begin{theorem}\label{lm:1}
If $f\in{S_{\alpha}}$ then $g(z):=z^{-\alpha}f(1/z)$ also belongs
to $S_{\alpha}$ and their representing measures are related by
$\mu_{g}=N_{\alpha}\mu_{f}$.
\end{theorem}
\textbf{Proof.} Indeed, using definition (\ref{eq:rho-rep}) we
get:
\[
z^{-\alpha}f(1/z)=z^{-\alpha}\int\limits_{(0,\infty)}\frac{z^{\alpha}\rho_{\alpha}(dt)}{(z+t)^{\alpha}}+
z^{-\alpha}\rho_0+\rho_{\infty}=\int\limits_{(0,\infty)}\frac{\rho_{\alpha}(dt)}{(z+t)^{\alpha}}+
z^{-\alpha}\rho_0+\rho_{\infty}
\]
which is precisely the representation (\ref{eq:mu-rep}).  Since
$\rho_{\alpha}\in\M_{\alpha}$ is arbitrary, the claim
follows.~$\square$

The next important result is due to Sokal \cite{Sokal}:
\begin{theorem}\emph{(Sokal, \cite{Sokal})}\label{lm:2}
A function $f$ defined on $(0,\infty)$ has holomorphic extension
$f\in{S_{\alpha}}$ if and only if
\begin{equation}\label{eq:Fnk}
F_{n,k}^{\alpha}(x):=(-1)^n D^{k}(x^{n+k+\alpha-1}D^nf(x))\geq{0}
\end{equation}
for all integers $n,k\geq{0}$ and all $x>0$.  Here $D=d/dx$.
\end{theorem}
\rem Differentiating (\ref{eq:mu-rep}) under the integral sign or
writing
\begin{equation}\label{eq:Fnkderivative}
(-1)^{n} D^{k}(x^{n+k+\alpha}D^n(-f'(x))=(-1)^{(n+1)}
D^{k}(x^{(n+1)+k+\alpha-1}D^{(n+1)}f(x)),
\end{equation}
we see that $f\in S_{\alpha}$ implies $-f'\in S_{\alpha+1}$, so
 that $(-1)^nf^{(n)}\in{S_{\infty}}$ for all integers
$n\geq{0}$ once $f\in{S_{\infty}}$.

\begin{theorem}\label{lm:3}
If $\alpha<\beta$  then $S_{\alpha}\subset{S_{\beta}}$. Moreover,
each $f\in{S_\alpha}$ defined by \emph{(\ref{eq:mu-rep})} can be
written as
\begin{equation}\label{eq:fviamubeta}
f(z)=\int\limits_{(0,\infty)}\frac{\mu_{\beta}(dy)}{(y+z)^{\beta}}+\mu_{\infty},
\end{equation}
where  $\mu_{\beta}\in\M_{\beta}$ and
\begin{equation}\label{eq:mualph-beta}
\mu_{\beta}(dy)=\frac{\Gamma(\beta)dy}{\Gamma(\alpha)\Gamma(\beta-\alpha)}\int\limits_{[0,y)}\frac{\mu_{\alpha}(du)}{(y-u)^{\alpha+1-\beta}}.
\end{equation}
Conversely, given $\mu_{\beta}\in\M_{\beta}$ which represents
$f\in{S_{\alpha}}$ we can recover $\mu_{\alpha}$ from
\begin{equation}\label{eq:mubeta-alpha}
F_{\mu_{\alpha}}(y)=\frac{\Gamma(\alpha)}{\Gamma(\beta)\Gamma(\alpha-\beta+n+1)}\left(\frac{d}{dy}\right)^n
\int\limits_{[0,y)}\frac{\mu_{\beta}(du)}{(y-u)^{\beta-\alpha-n}}
\end{equation}
with $n=[\beta-\alpha]$ and
$\mu_{\alpha}(\{\infty\})=\mu_{\beta}(\{\infty\})$.
\end{theorem}
\rem Formula (\ref{eq:mualph-beta}) generalizes \cite[Chapter
VIII, Corollary 3a.1, p.330]{Widder} which in our notation
connects the measures $\mu_{1}$ and $\mu_{2}$.  The transformation
$\mu_{\alpha}\to\mu_{\beta}$ defined in (\ref{eq:mualph-beta}) is
the left-sided Riemann-Liouville fractional integral.  Its precise
definition and inversion are investigated in the Appendix.

\textbf{Proof.} Since
\begin{equation}\label{eq:uzalpha}
\frac{1}{(u+z)^{\alpha}}=\frac{\Gamma(\beta)}{\Gamma(\alpha)\Gamma(\beta-\alpha)}
\int\limits_{(0,\infty)}\frac{t^{\beta-\alpha-1}dt}{(z+u+t)^{\beta}},
\end{equation}
we have
\begin{multline}\label{eq:alpha-beta-chain1}
\int\limits_{[0,\infty)}\frac{\mu_{\alpha}(du)}{(u+z)^{\alpha}}
=\frac{\Gamma(\beta)}{\Gamma(\alpha)\Gamma(\beta-\alpha)}\int\limits_{[0,\infty)}\mu_{\alpha}(du)
\int\limits_{(0,\infty)}\frac{t^{\beta-\alpha-1}dt}{(z+u+t)^{\beta}}\biggl|_{y=u+t}
\\
=\frac{\Gamma(\beta)}{\Gamma(\alpha)\Gamma(\beta-\alpha)}\int\limits_{[0,\infty)}\mu_{\alpha}(du)
\int\limits_{(u,\infty)}\frac{(y-u)^{\beta-\alpha-1}dy}{(z+y)^{\beta}}
\\
=\frac{\Gamma(\beta)}{\Gamma(\alpha)\Gamma(\beta-\alpha)}
\int\limits_{(0,\infty)}\frac{dy}{(z+y)^{\beta}}\int\limits_{[0,y)}(y-u)^{\beta-\alpha-1}\mu_{\alpha}(du)
=\int\limits_{[0,\infty)}\frac{\mu_{\beta}(dy)}{(z+y)^{\beta}},
\end{multline}
where $\mu_{\beta}(dy)$ is defined by (\ref{eq:mualph-beta}) and
we can use Tonelli's theorem \cite[Chapter~20, Corollary~7]{RF} or
\cite[Theorem~3.4.5]{Bogachev} to justify the interchange of
integrations. Also setting $z=1$ we see by Tonelli's theorem and
condition (\ref{eq:mu-condition}) that the measure $\mu_{\beta}$
belongs to  $\M_\beta$.  According to Remark~A1 in the Appendix
the measure $\mu_{\beta}$ has no atom at zero.  This allows us to
remove zero from the domain of integration in
(\ref{eq:fviamubeta}). This proves (\ref{eq:fviamubeta}) and the
inclusion $S_{\alpha}\subset{S_{\beta}}$. The inversion formula
(\ref{eq:mubeta-alpha}) is Theorem~\ref{th:Leftsided} in the
Appendix.~~$\square$

\begin{theorem}\label{lm:31}
Each $f\in{S_\alpha}$ defined by \emph{(\ref{eq:rho-rep})} can
also be written as
\[
f(z)=\int\limits_{(0,\infty)}\frac{\rho_{\beta}(dx)}{(1+xz)^{\beta}}+\rho_0,
\]
where  $\rho_{\beta}\in\M_{\beta}$ and
\begin{equation}\label{eq:rhoalph-beta}
\rho_{\beta}(dx)=\frac{\Gamma(\beta)x^{\alpha-1}dx}{\Gamma(\alpha)\Gamma(\beta-\alpha)}
\biggl\{\int\limits_{(x,\infty)}\frac{u^{1-\beta}\rho_{\alpha}(du)}{(u-x)^{\alpha-\beta+1}}+\rho_{\infty}\biggr\}.
\end{equation}
Conversely, given $\rho_{\beta}\in\M_{\beta}$ which represents
$f\in{S_{\alpha}}$ we can recover $\rho_{\alpha}$ from
\begin{multline}\label{eq:rhobeta-alph}
F_{\rho_{\alpha}}(y)=\rho_0+\frac{\Gamma(\alpha)}{\Gamma(\beta)\Gamma(\alpha-\beta+n+1)}\times
\\
\biggl\{\!\alpha\!\!\int\limits_{(0,y)}\!\!x^{\alpha-1}dx
\left(-x^2\frac{d}{dx}\right)^nx^{\beta-\alpha-n}\!\!\!\!\int\limits_{(x,\infty)}\!\!\!\!\frac{\rho_{\beta}(ds)}{s^{\alpha+n}(s-x)^{\beta-\alpha-n}}
-y^{\alpha}\left(-y^2\frac{d}{dy}\right)^ny^{\beta-\alpha-n}\!\!\!\!\int\limits_{(y,\infty)}\!\!\!\!\frac{\rho_{\beta}(ds)}{s^{\alpha+n}(s-y)^{\beta-\alpha-n}}\!\biggr\}
\end{multline}
with $n=[\beta-\alpha]$ and
$\rho_{\alpha}(\{\infty\})=[\Gamma(\beta-\alpha)\Gamma(\alpha+1)/\Gamma(\beta)]\lim\limits_{y\to\infty}y^{-\alpha}F_{\rho_{\beta}}(y)$.
\end{theorem}
\rem  The transformation $\rho_{\alpha}\to\rho_{\beta}$ defined in
(\ref{eq:rhoalph-beta}) is the right-sided Kober-Erdelyi
fractional integral. Its precise definition and inversion are
investigated in the Appendix.

\textbf{Proof.} To demonstrate (\ref{eq:rhoalph-beta}) we employ
the connection formula
$\rho_{\beta}^*(dy)=y^{-\beta}\mu_{\beta}(dy)$, where
$\rho_{\beta}^*$ is the image of $\rho_{\beta}$ under
$y\to{y^{-1}}$. We have
($A=\Gamma(\beta)/[\Gamma(\alpha)\Gamma(\beta-\alpha)]$ for
brevity)
\begin{multline*}
y^{\beta}\rho_{\beta}^*(dy)=\mu_{\beta}(dy)=Ady\int\limits_{(0,y)}\frac{\mu_{\alpha}(du)}{(y-u)^{\alpha+1-\beta}}
+\mu_{0}Ay^{\beta-\alpha-1}dy
\\
=Ady\!\!\!\int\limits_{(1/y,\infty)}\!\!\!\frac{t^{\alpha+1-\beta}\mu_{\alpha}^{*}(dt)}{(yt-1)^{\alpha+1-\beta}}+A\rho_{\infty}y^{\beta-\alpha-1}dy
=Ay^{\beta-\alpha-1}dy\!\!\!\int\limits_{(1/y,\infty)}\!\!\!\frac{t^{1-\beta}\rho_{\alpha}(dt)}{(t-1/y)^{\alpha+1-\beta}}+A\rho_{\infty}y^{\beta-\alpha-1}dy.
\end{multline*}
Dividing by $y^{\beta}$ and making substitution $y=1/x$,
$dy=-dx/x^2$ we arrive at (\ref{eq:rhoalph-beta}).  The inversion
formula (\ref{eq:rhobeta-alph}) is Theorem~\ref{th:Rightsided} in
the Appendix.~~$\square$

The following result is also due to Sokal \cite[formulas (10a) and
(10b)]{Sokal}.
\begin{theorem}\label{lm:4}
$\bigcap_{\alpha>0}S_{\alpha}=\{\text{\rm non-negative
constants}\}$.
\end{theorem}
This result suggests the following definition:
$$
S_0:=\{\text{\rm non-negative constants}\}.
$$

Recall that a function $f:~(0, \infty)\to\mathbb{R}$ is said to be
completely monotonic if $f$ has derivatives of all orders and
satisfies $(-1)^nf^{(n)}(x)\ge{0}$, for all $x>0$ and
$n=0,1,\ldots$ We denote the set of completely monotonic functions
by ${\cal CM}$. The following result was pointed out to us by
Christian Berg and is also hinted at in \cite{Sokal}.
\begin{theorem}\label{lm:union}
$\overline{S_{\infty}}={\cal CM}$, where the closure is taken with
respect to pointwise convergence on $(0,\infty)$.
\end{theorem}
\textbf{Proof.} The inclusion $\overline{S_{\infty}}\subset{\cal
CM}$ follows from the fact that each
$f\in\bigcup_{\alpha>0}S_{\alpha}$ is completely monotonic
combined with closedness of  ${\cal CM}$ under pointwise
convergence \cite[Corollary~1.6]{SSV}.  To prove the reverse
inclusion recall that according to Bernstein's theorem (see
\cite[Chapter IV, Theorem~12b]{Widder} or \cite[Theorem~1.4]{SSV})
each completely monotonic function is the Laplace transform of a
nonnegative measure:
$$
f(x)=\Laplace(\sigma;x):=\int\limits_{[0,\infty)}e^{-xt}\sigma(dt).
$$
Since $f(x)$ is non-increasing the sequence $f(1/n)$ is
non-decreasing.  Hence, two cases are possible: 1) $f(1/n)$ is
bounded by a constant $Ñ$ or  2)
$\lim_{n\to\infty}f(1/n)=+\infty$. Define $a_n:=\sqrt{n}$ in the
first case and  $a_n:=f(1/n)$ in the second. Define the following
sequence of functions:
$$
f_{n}(x)=\int\limits_{[0,\infty)}\left(1+\frac{xt}{a_n^2}\right)^{-a_n^2}e^{-\frac{t}{n}}{\sigma(dt)}.
$$
Obviously,
$$
f_n(x)\leq\int\limits_{[0,\infty)}e^{-\frac{t}{n}}{\sigma(dt)}=f(1/n),
$$
so that the integral defining $f_n(x)$ exists. Moreover,
$$
f_n(x)=\int\limits_{[0,\infty)}\frac{\sigma_n(du)}{(1+xu)^{a_n^2}},
$$
where $\sigma_n(du)=a_n^2e^{\frac{-a_n^2u}{n}}\sigma(du)$.
Clearly, $f_n\in\bigcup_{\alpha>0}S_{\alpha}$. Next, we have
$$
f_n(x)-f(x)=\int\limits_{[0,\infty)}\left(\left(1+\frac{xt}{a_n^2}\right)^{-a_n^2}-e^{-xt}+e^{-xt}-e^{-(x-1/n)t}\right)e^{-\frac{t}{n}}{\sigma(dt)}
$$
\begin{equation}\label{eq:diifference}
f_n(x)-f(x)=\int\limits_{[0,\infty)}\left(\left(1+\frac{xt}{a_n^2}\right)^{-a_n^2}-e^{-xt}\right)e^{-\frac{t}{n}}{\sigma(dt)}+f(x+1/n)-f(x).
\end{equation}
It is verified by straightforward calculus that the maximum of
$$
\psi(t):=\left|\left(1+\frac{xt}{a_n^2}\right)^{-a_n^2}-e^{-xt}\right|
$$
is attained at the point $t^*$ satisfying
$$
\left(1+\frac{xt}{a_n^2}\right)^{-a_n^2-1}=e^{-xt}.
$$
Hence,
$$
\max_t
|\psi(t)|=\frac{xt^*e^{-xt^*}}{a_n^2}\leq\frac{e^{-1}}{a_n^2}.
$$
It follows from (\ref{eq:diifference}) that
$$
|f_n(x)-f(x)|\leq\frac{e^{-1}}{a_n^2}f(1/n)+|f(x+1/n)-f(x)|.
$$
Due to the definition of $a_n$ and continuity of $f(x)$ we
conclude that $\lim_{n\to\infty}f_n(x)=f(x)$ for each
$x>0$.~~$\square$

\begin{theorem}\label{lm:6}
If  $f\in{S_{\alpha}}$ and $g\in{S_{\beta}}$ then
$fg\in{S_{\alpha+\beta}}$.
\end{theorem}
\textbf{Proof.} See \cite[Chapter VII, paragraph 7.4]{HW}.

\rem Theorems~\ref{lm:3} and \ref{lm:6} show that the union
$S_{\infty}$ is a cone with multiplication: if
$f,g\in{S_{\infty}}$ then $af+bg\in{S_{\infty}}$, for all
$a,b\geq{0}$ and $fg\in{S_{\infty}}$.

\begin{theorem}\label{lm:7}
Each
\[
f(z)=\int\limits_{[0,\infty)}\frac{\mu_{\alpha}(du)}{(u+z)^{\alpha}}+\mu_{\infty}\in{S_{\alpha}}
\]
can be represented in the form
\begin{equation}\label{eq:flaplace}
f(z)=\frac{1}{\Gamma(\alpha)}\Laplace(u^{\alpha-1}\Laplace(\mu_{\alpha};u)du;z)+\mu_{\infty},
\end{equation}
 where $\Laplace$ denotes the Laplace transform.
\end{theorem}
\textbf{Proof.}  Write
\[
\frac{1}{(u+z)^{\alpha}}=\frac{1}{\Gamma(\alpha)}\int\limits_{[0,\infty)}e^{-(u+z)t}t^{\alpha-1}dt
\]
and apply Tonelli's theorem to show that the iterated integral in
(\ref{eq:flaplace}) exists and is equal to $f(z)$.~~$\square$

\noindent\rem Formula (\ref{eq:flaplace}) has been found in
\cite{Yur} for absolutely continuous measures.

Introduce the standard notation for the right-sided
Riemann-Liouville fractional integral \cite[section~2.2]{KST},
\cite[\S5]{SKM} and the right-sided Caputo fractional derivative
\cite[section~2.4]{KST}, \cite[section 2.4.1]{Podlubny}:
\[
I_{\lambda}^{-}f:=\frac{1}{\Gamma(\lambda)}\int\limits_{x}^{\infty}\frac{f(t)dt}{(t-x)^{1-\lambda}},
\]
\[
{}^C\!D_{\lambda}^{-}f=(I_{\lambda}^{-})^{-1}f=\frac{(-1)^n}{\Gamma(n-\lambda)}
\int\limits_{x}^{\infty}\frac{f^{(n)}(t)dt}{(t-x)^{1+\lambda-n}},~~n=[\lambda]+1.
\]

\begin{theorem}\label{lm:8}
For a fixed non-negative measure $\mu\in\M_{\alpha}$  and
$\beta>\alpha>0$ denote
\begin{equation}\label{eq:falpha}
f_{\alpha}(z)=\int\limits_{[0,\infty)}\frac{\mu(du)}{(u+z)^{\alpha}}+\mu_{\infty},
\end{equation}
and
\begin{equation}\label{eq:fbeta}
f_{\beta}(z)=\int\limits_{[0,\infty)}\frac{\mu(du)}{(u+z)^{\beta}}+\mu_{\infty}.
\end{equation}
Then for all $x>0$
\begin{equation}\label{eq:falphabeta}
f_{\alpha}(x)
=\frac{\Gamma(\beta)}{\Gamma(\alpha)\Gamma(\beta-\alpha)}\int\limits_{x}^{\infty}\frac{(f_{\beta}(t)-\mu_{\infty})dt}{(t-x)^{1+\alpha-\beta}}+\mu_{\infty}
=\frac{\Gamma(\beta)}{\Gamma(\alpha)}I_{\beta-\alpha}^{-}[f_{\beta}(t)-\mu_{\infty}](x)+\mu_{\infty}
\end{equation}
and
\begin{equation}\label{eq:fbetaalpha}
f_{\beta}(x)
=\frac{\Gamma(\alpha)(-1)^n}{\Gamma(\beta)\Gamma(n-\beta+\alpha)}\int\limits_{x}^{\infty}\frac{f_{\alpha}^{(n)}(t)dt}{(t-x)^{1+\beta-\alpha-n}}
+\mu_{\infty}=\frac{\Gamma(\alpha)}{\Gamma(\beta)}{}^C\!D_{\beta-\alpha}^{-}[f_{\alpha}(t)](x)+\mu_{\infty},
\end{equation}
 where $n=[\beta-\alpha]+1$.
\end{theorem}
\textbf{Proof.} To prove (\ref{eq:falphabeta}) substitute
(\ref{eq:fbeta}) for $f_{\beta}$ into (\ref{eq:falphabeta}) and
exchange the order of integration which is legitimate by Tonelli's
theorem again. Conditions $\beta>\alpha>0$ guarantee the existence
of the inner integral.  Formula (\ref{eq:fbetaalpha}) is one of
several inversion formulas for the Riemann-Liouville integral in
(\ref{eq:falphabeta}) which is applicable since $f_{\alpha}$ is
infinitely differentiable   (see, for instance,
\cite[section~2.4]{KST}, \cite[section 2.4.1]{Podlubny}). To
demonstrate its validity differentiate under integral sign
($n\geq{1}$)
$$
f^{(n)}_{\alpha}(t)=(-1)^n(\alpha)_{n}\int\limits_{[0,\infty)}\frac{\mu(du)}{(u+t)^{\alpha+n}},
$$
where
$(\alpha)_{n}=\alpha(\alpha+1)\cdots(\alpha+n-1)=\Gamma(\alpha+n)/\Gamma(\alpha)$,
substitute into (\ref{eq:fbetaalpha}) and exchange the order of
integrations.~~$\square$

\rem The Riemann-Liouville fractional derivative cannot be used in
(\ref{eq:fbetaalpha}) since, in general, the resulting integral
would diverge.

\begin{theorem}\label{lm:pointwise}
The class $S_{\alpha}$, $\alpha>1$ is closed under pointwise
limits: if $\{f_n\}_{n=1}^{\infty}\subset{S_{\alpha}}$ and if the
limit $\lim\limits_{n\to\infty}f_n(x)=f(x)$ exists for all $x>0$
then $f\in{S_{\alpha}}$.
\end{theorem}
\textbf{Proof.} A proof for $\alpha=1$ can be found in
\cite[Proposition~1]{Hirsch} or by a different argument in
\cite[Theorem~2.2(iii)]{SSV}.  The latter carries over
\emph{mutatis mutandis} to all $\alpha>0$. ~~$\square$

\bigskip

\paragraph{3. Complex variable properties.}
Clearly, $f$ is holomorphic in $\C\!\setminus\!(-\infty,-r]$,
where $r=\inf\{x:~x\in\mathbf{supp}(\mu_{\alpha})\}$ and
$f(\overline{z})=\overline{f(z)}$. In particular, if
$0\notin\mathbf{supp}(\mu_{\alpha})$ (or, equivalently,
$\mathbf{supp}(\rho_{\alpha})$ is bounded) then the function $f$
can be represented by the power series
\[
f(z)=\sum\limits_{k=0}^{\infty}(-1)^k\frac{(\alpha)_k}{k!}\rho_k(\alpha)z^k,
\]
convergent in the disk $|z|<1/R$. Here
$R=\sup\{x:~x\in\mathbf{supp}(\rho_{\alpha})\}$,
$(\alpha)_k=\Gamma(\alpha+k)/\Gamma(\alpha)$, and
\begin{equation}\label{eq:rho-finite}
\rho_k(\alpha)=\int\limits_{[0,R]}t^kd\rho_{\alpha}(t)<\infty,~~k=0,1,\ldots
\end{equation}
are the moments of the measure $\rho_{\alpha}$ which are finite
due to (\ref{eq:rho-condition}).

For functions belonging to $S_1$ Krein \cite[Appendix,
Theorem~4]{KN} found the following celebrated characterization.
\begin{theorem}\label{th:Krein}
A function $f$ holomorphic in the cut plane
$\C\setminus\!(-\infty,0]$ belongs to $S_1$ iff $f(x)\geq{0}$ for
$x>0$ and $\Im{f(z)}\leq{0}$ for $\Im{z}>0$.
\end{theorem}

Because of the special role played by $S_1$ (including a number of
stability properties \cite{Berg}, connection to  continued
fractions \cite{CPVWJ} and Pad\'{e} approximation \cite{BG}) it is
interesting to relate the functions from $S_{\alpha}$ to $S_1$.
One way of doing this is provided by Theorem~\ref{lm:8} (it works
if $\mu\in\M_{1}$), another approach is presented in the following
two theorems.

\begin{theorem}\label{th:Salpha1}
Suppose $f\in{S_{\alpha}}$, $0<\alpha\leq{1}$. Then
$f^{1/\alpha}\in{S_1}$.
\end{theorem}
\textbf{Proof.} For $\Im{z}>0$ and $t>0$ we have ($\arg(z)$
denotes the principal value of the argument of $z$):
\[
-\pi\alpha<\arg(z+t)^{-\alpha}<0~~\Rightarrow~~
-\pi\alpha<\arg(f(z))\leq{0}~~\Rightarrow~~\Im(f(z)^{1/\alpha})\leq{0}.
\]
Since $f^{1/\alpha}(z)$ is holomorphic in $\C\setminus(-\infty,0]$
and non-negative for $z>0$ we get the conclusion by Krein's
theorem~\ref{th:Krein}.~~$\square$

\begin{theorem}\label{th:Salpha2}
Suppose $f\in{S_{\alpha}}$, $\alpha\geq{1}$. Then
$g(z):=f(z^{1/\alpha})\in{S_1}$.
\end{theorem}
\textbf{Proof.}  Indeed $g(z)$ is  holomorphic in
$\C\setminus(-\infty,0]$ and non-negative for $x>0$.  Next for
$\Im{z}>0$ and $t>0$ we have ($\arg(z)$ denotes the principal
value of the argument of $z$):
\[
0<\arg(z^{1/\alpha})<\pi/\alpha~~\Rightarrow~~0<\arg(z^{1/\alpha}+t)<\pi/\alpha,
~~\Rightarrow~~\Im(z^{1/\alpha}+t)^{-\alpha}<0
\]
Integrating the last expression with respect to non-negative
measure preserves the lower half plane, so that the proof is
completed by Krein's theorem~\ref{th:Krein}.~~$\square$

\rem For $\alpha>1$ the mapping $S_{\alpha}\to{S_1}$ defined by
$f(z)\to{g(z)}:={f(z^{1/\alpha})}$ is clearly not surjective as
can be seen immediately by taking $g(z)=1/(1+z)\in{S_1}$,
$f(z)=g(z^{\alpha})\notin{S_{\alpha}}$ (since it is not
holomorphic in the upper half-plane).  The conditions (a) $f(z)$
is holomorphic in the cut plane $\C\setminus(-\infty,0]$, (b)
$f(x)\geq{0}$ for $x>0$ and (c) $\Im{f(z)}\leq{0}$ for
$0<\arg(z)<\pi/\alpha$ are necessary for $f$ to belong to
$S_{\alpha}$.  Unfortunately, these conditions are not sufficient
as the following example shows:
$$
f(z)=\frac{1}{(z+1)^2}-\frac{1}{2(z+2)^2}.
$$
Indeed, a straightforward computation yields ($z=x+iy$):
$$
\Im{f(z)}=-\frac{y(30+87x+96x^2+50x^3+12x^4+x^5+12y^2+22xy^2+12x^2y^2+2x^3y^2+xy^4)}{(1+2x+x^2+y^2)^2(4+4x+x^2+y^2)^2}.
$$
Hence, for $\{z\in\C:~\Re{z}>0,\Im{z}>0\}$ we get $\Im{f(z)}<0$
and $f(x)\geq{0}$ for $x>0$. However, $f\notin{S_2}$ since the
representing measure is signed.  In fact, $f\in{S_3}$ since
$$
f(z)=\int\limits_{1}^{\infty}\frac{d\nu(t)}{(z+t)^3},~\text{where}~d\nu(t)\!=\!\left\{\begin{array}{l}\!\!2dt,~t\in(1,2],\\\!\!dt,~t\in(2,\infty).\end{array}\right.
$$
We thank Alex Gomilko for this example. It provides a partial
answer to the question asked by Sokal at the end of \cite{Sokal}.

\begin{theorem}
Each $f\in{S_{\alpha}}$ satisfies
\begin{equation}\label{eq:falpha-estimate}
|f(z)|\leq A\left|\frac{z-1}{\Im(z)}\right|^{\alpha}+\mu_{\infty},
\end{equation}
where $A=\int_{[0,\infty)}\mu_{\alpha}(du)/(u+1)^{\alpha}<\infty$,
and $\Re(z)\leq{0}$.
\end{theorem}
\textbf{Proof.} We have
\begin{multline*}
\left|f(z)\right|=\left|\,\int\limits_{[0,\infty)}\frac{\mu_{\alpha}(du)}{(u+z)^{\alpha}}+\mu_{\infty}\right|
=\left|\,\int\limits_{[0,\infty)}\frac{\mu_{\alpha}(du)}{(u+1)^{\alpha}}\left(\frac{u+1}{u+z}\right)^{\alpha}+\mu_{\infty}\right|
\\
\leq
A\max\limits_{u\geq{0}}\left|\frac{u+1}{u+z}\right|^{\alpha}+\mu_{\infty}=A[\psi(z)]^{\alpha/2}+\mu_{\infty},
\end{multline*}
with $A$ given above and
\[
\psi(z)=\max\limits_{u\geq{0}}\left|\frac{u+1}{u+z}\right|^2=\left|\frac{z-1}{\Im(z)}\right|^2.
\]
The last equality is true for $\Re(z)\leq{0}$ by standard
calculus. ~~$\square$

\rem The estimate (\ref{eq:falpha-estimate}) may seem very weak
because the right hand side does not tend to zero as
$|z|\to\infty$ while $f(z)\to{0}$. However, the decrease of $f$
may be arbitrarily slow so that it is difficult to expect a better
estimate valid for the whole class $S_{\alpha}$.

\bigskip

\paragraph{4. Exact Stieltjes order.} We will say that $f$ is of exact
Stieltjes order $\alpha^*$ if $f\in{S_{\infty}}$ and
\begin{equation}\label{eq:exactorder}
\alpha^*[f]=\inf\{\alpha:~f\in{S_{\alpha}}\}.
\end{equation}

\begin{theorem}\label{th:order}
If  $f$ is of exact Stieltjes order $\alpha^*$ then
$f\in{S_{\alpha^*}}$. In other words the infimum  in
\emph{(\ref{eq:exactorder})} is always attained.
\end{theorem}
\textbf{Proof.} Since $f\in{S_{\infty}}$ it is infinitely
differentiable  so that $F_{n,k}^{\alpha}(x)$ in (\ref{eq:Fnk}) is
well defined and continuous in $\alpha$ for each fixed $x>0$ for
all $\alpha>\alpha^*$.  Passing to the limit $\alpha\to\alpha^*$
we verify that (\ref{eq:Fnk}) is true for $\alpha=\alpha^*$.
Hence, by Sokal's theorem~\ref{lm:2}
$f\in{S_{\alpha^*}}$.~~$\square$

\begin{theorem}\label{th:criterium}
Suppose $f\in{S_{\beta}}$.  Then $\beta>\alpha^*[f]$  iff the
function
\begin{equation}\label{eq:criterium}
\Phi(y)=\int\limits_{(0,y)}\frac{\mu_{\beta}(du)}{(y-u)^{\varepsilon}}
\end{equation}
is non-decreasing on $(0,\infty)$  for some
$\varepsilon\in(0,\min\{\beta,1\})$.
\end{theorem}

\textbf{Proof.} Assume $\Phi(y)$ is non-decreasing. We need to
show that $\beta$ is not exact. According to Definition~1 in the
Appendix and Theorem~\ref{th:Leftsided}
$$
I^{+}_{1-\varepsilon}\mu_{\beta}=\frac{\Phi(y)dy}{\Gamma(1-\varepsilon)}+\mu_{\infty}\delta_{\infty}\in\M_{\beta+1-\varepsilon}=\M_{\alpha+1},
$$
where $\alpha=\beta-\varepsilon$ and
$\mu_{\infty}=\mu_{\beta}(\{\infty\})$. Hence,
$$
\int\limits_{[0,\infty)}\frac{\Phi(y)dy}{(1+y)^{\alpha+1}}<\infty.
$$
This implies that
\begin{equation}\label{eq:dphi-finite}
\int\limits_{[0,\infty)}\frac{d\Phi(y)}{(1+y)^{\alpha}}<\infty.
\end{equation}
Indeed, for each $t>0$  integration by parts,
$$
\frac{\Phi(y)}{\alpha(1+y)^{\alpha}}\biggl|^{t}_{0}=\frac{1}{\alpha}\int\limits_{[0,t)}\frac{d\Phi(y)}{(1+y)^{\alpha}}-\int\limits_{[0,t)}\frac{\Phi(y)dy}{(1+y)^{\alpha+1}},
$$
shows that $\lim_{t\to\infty}\Phi(t)(1+t)^{-\alpha}$ exists
(finite or infinite). This limit must be zero since otherwise
$\Phi(y)(1+y)^{-\alpha}>C>0$ for all $y>M$ and
$$
\int\limits_{[0,\infty)}\frac{\Phi(y)dy}{(1+y)^{\alpha+1}}>C\int\limits_{[M,\infty)}\frac{dy}{1+y}=\infty.
$$
Hence, $\lim_{t\to\infty}\Phi(t)(1+t)^{-\alpha}=0$ and the above
integration by parts  proves (\ref{eq:dphi-finite}).
 It follows that the measure $\mu_{\alpha}$  whose
distribution function is equal to $A\Phi(y)$ and whose atom at
infinity is equal to $\mu_{\infty}$ belongs to $\M_{\alpha}$. Here
$A=\Gamma(\alpha)/[\Gamma(\beta)\Gamma(\alpha-\beta+1)]$. Consider
the function
\[
g(z)=\int\limits_{[0,\infty)}\frac{\mu_{\alpha}(du)}{(u+z)^{\alpha}}+\mu_{\infty}.
\]
By Theorem~\ref{lm:3} we have $g\in{S_{\beta}}$ and
\[
g(z)=\int\limits_{[0,\infty)}\frac{\tilde{\mu}_{\beta}(du)}{(u+z)^{\beta}}+\mu_{\infty},
\]
where $\tilde{\mu}_{\beta}(du)$ is given by
(\ref{eq:mualph-beta}). But then $\mu_{\alpha}$ and
$\tilde{\mu}_{\beta}$ are related by (\ref{eq:mubeta-alpha}) which
coincides with (\ref{eq:criterium}) times $A$ (note that $n=0$ in
(\ref{eq:mubeta-alpha}) because $\beta-\alpha=\varepsilon<1$).
This proves that $\tilde{\mu}_{\beta}=\mu_{\beta}$ so that
$f(z)=g(z)\in{S_{\alpha}}$, $\alpha<\beta$.

Conversely, if $\beta$ is not exact choose
$\varepsilon\in(0,\beta-\alpha^*)$, $\varepsilon<1$. We have
$f\in{S_{\beta-\varepsilon}}$. According to Theorem~\ref{lm:3} the
function $A\Phi(y)$, where $\Phi$ is defined in
(\ref{eq:criterium}) equals the distribution function of the
measure $\mu_{\beta-\varepsilon}$ and so is non-decreasing.
~$\square$

\begin{corollary}\label{cr:limit}
Suppose $f\in{S_{\beta}}$, $\varepsilon\in(0,\min\{\beta,1\})$ and
the following limit exists:
\[
\lim\limits_{y\to+\infty}\frac{\Phi(2y)}{\Phi(y)}=A,
\]
where $\Phi$ is defined by \emph{(\ref{eq:criterium})}. If $A<1$
then $\beta$ is the exact Stieltjes order of $f$.
\end{corollary}
\textbf{Proof.} Clearly, the condition $A<1$ implies that
$\Phi(y)$ cannot be non-decreasing so that by
Theorem~\ref{th:criterium} $\beta$ must be exact.~$\square$

\begin{corollary}\label{cr:compact}
Suppose $f\in{S_{\beta}}$ and the support of the measure
$\mu_{\beta}$ is compact.  Then $\beta$ is the exact Stieltjes
order of $f$.
\end{corollary}
\textbf{Proof.}  Indeed, for all $y>B$,
$B:=\sup\{x:~x\in\mathbf{supp}(d\mu)\}$, the function $\Phi(y)$ is
strictly decreasing for each $\varepsilon\in(0,\min\{\beta,1\})$,
so that by Theorem~\ref{th:criterium} $\beta$ must be
exact.~$\square$

Consider three prototypical examples.

\textbf{Example~1}. Find the exact Stieltjes order of ($\alpha>1$)
\[
f(z)=\int\limits_{0}^{1}\frac{dt}{(z+t)^{\alpha}}
=\frac{1}{\alpha-1}\left(\frac{1}{z^{\alpha-1}}-\frac{1}{(1+z)^{\alpha-1}}\right).
\]
Method I: by corollary~\ref{cr:compact} $\alpha^*=\alpha$ since
$\mathbf{supp}(d\mu)$ is compact.

Method II: by theorem~\ref{th:criterium} compute $\Phi(y)$. Let
$I(A)$ be the indicator function of a set $A$. We have
($0<\varepsilon<1$)
\[
\Phi(y)=\int\limits_{0}^{y}\frac{I([0,1])du}{(y-u)^{\varepsilon}}=
\left\{\!\!\!
\begin{array}{l}
y^{1-\varepsilon}/(1-\varepsilon),~~~0<y\leq{1},\\[5pt]
[y^{1-\varepsilon}-(y-1)^{1-\varepsilon}]/(1-\varepsilon),~~~y>1.
\end{array}
\right.
\]
It is straightforward to check that this function is decreasing
for $y>1$ so that again $\alpha^*=\alpha$.

\textbf{Example~2.} Find the exact Stieltjes order of ($\alpha>1$)
\[
f(z)=\int\limits_{1}^{\infty}\frac{dt}{(z+t)^{\alpha}}.
\]
By theorem~\ref{th:criterium} compute
\[
\Phi(t)=\int\limits_{0}^{t}\frac{I([1,\infty))du}{(t-u)^{\varepsilon}}=
\left\{\!\!\!
\begin{array}{l}
0,~~~0<t\leq{1} \\[5pt]
(t-1)^{1-\varepsilon}/(1-\varepsilon),~~~t>1
\end{array}
\right.
\]
This function is non-decreasing on $[0,\infty)$, so that
$\alpha^*<\alpha$.  To find the exact order we compute
\[
f(z)=\frac{1}{(\alpha-1)(1+z)^{\alpha-1}}=\int\limits_{0}^{\infty}\frac{d\nu(t)}{(z+t)^{\alpha-1}},
\]
where the measure $d\nu(t)$ is concentrated at one point $t=1$
with $\nu(\{1\})=(\alpha-1)^{-1}$ so that by
Corollary~\ref{cr:compact} $\alpha^*=\alpha-1$.

\textbf{Example~3.} According to Euler's integral representation
\cite[Theorem~2.2.1]{AAR} the Gauss hypergeometric functions
$_2F_1$ can be written as
$$
f(z):=_2F_1(a,b;c;-z)=\frac{\Gamma(c)}{\Gamma(b)\Gamma(c-b)}\int_0^1\frac{u^{b-1}(1-u)^{c-b-1}}{(1+zu)^a}du,
~~c>b>0.
$$
Assume that $0<a\leq{b}$. By the above formula $f\in{S_{a}}$.
Change of variable $u=1/t$ yields
$$
_2F_1(a,b;c;-z)=\frac{\Gamma(c)}{\Gamma(b)\Gamma(c-b)}\int_1^\infty\frac{\mu(t)dt}{(z+t)^a},
$$
where
$$
\mu(t)={t^{a-c}}{(t-1)^{c-b-1}},~~~~t>1.
$$
We aim to show that $a$ is the exact Stieltjes order of $f$.
Computation gives
$$
\Phi(y):=\int\limits_{1}^{y}\frac{\mu(t)dt}{(y-t)^{\varepsilon}}=A(y-1)^{c-b-\varepsilon}{_2F_1}(c-a,c-b;c-b-\varepsilon+1;-(y-1)),
$$
where
$A=\Gamma(c-b)\Gamma(1-\varepsilon)/[\Gamma(c-b-\varepsilon+1)]$
for any given $\varepsilon\in(0,\min\{1,c-b\})$. Using  the
asymptotic formula \cite[formula (2.3.12)]{AAR} for ${_2F_1}$ as
$y\to\infty,$ we get
$$
\Phi(y)=A(y-1)^{-\varepsilon}\left(D+o(1)\right),
$$
where $D$ is a constant.  This implies that
$\lim_{y\to\infty}[\Phi(2y)/\Phi(y)]=2^{-\varepsilon}<1$ so that
by Corollary~\ref{cr:limit} the exact order is equal to $a$.

In our forthcoming paper \cite{KaPr} we will investigate the exact
Stieltjes order of the generalized hypergeometric function
$_{q+1}F_{q}$.

\rem  According to Theorem~\ref{lm:pointwise} the class
$S_{\alpha}$ is closed under pointwise limits. It is then
reasonable to ask whether the exact Stieltjes order is also
preserved by such limits.  The following example shows that the
answer is negative in general.  Consider the sequence of functions
($\alpha>1$)
$$
f_m(z)=\int\limits_{1}^{m}\frac{dt}{(z+t)^{\alpha}}=\frac{1}{(\alpha-1)(z+1)^{\alpha-1}}-\frac{1}{(\alpha-1)(z+m)^{\alpha-1}},~~~m=2,3,\ldots
$$
According to Corollary~\ref{cr:compact} each $f_m$ has exact order
$\alpha$ while pointwise
$$
\lim\limits_{m\to\infty}f_m(z)=\frac{1}{(\alpha-1)(z+1)^{\alpha-1}}~\text{for
each}~z\in\C\!\setminus\!(-\infty,0].
$$
According, to Example~2 above, the limit function is of exact
order $\alpha-1$.

\rem  Theorem~\ref{lm:8} shows that if $\alpha$ is the exact
Stieltjes order of $f$ then its fractional derivative of order
$\gamma$ will have the exact order $\alpha+\gamma$ while its
fractional integral of order $\gamma$ will have the exact order
$\alpha-\gamma$ provided that $\alpha>\gamma$ and
$\mu\in\M_{\alpha-\gamma}$.

\newpage

\begin{center}
\Large{\textbf{APPENDIX}}\\ \Large{\textbf{Inversion formulas for
some fractional integrals of measures on half-axis}}
\end{center}

\renewcommand{\thetheorem}{A\arabic{theorem}}
\setcounter{theorem}{0}
\renewcommand{\rem}{\par\refstepcounter{remark}\noindent\textbf{Remark
A\arabic{remark}.} }
\setcounter{remark}{0}

\bigskip

\bigskip

In this appendix we prove some facts about fractional integrals
and derivatives of Borel measures supported on $\R^+$.  Fractional
calculus is certainly a classical subject with a number of great
monographs available today, including
\cite{Diet,KST,Kiryakova,MR,Podlubny,SKM}. Fractional integrals
and derivatives have been studied in (weighted) spaces of
integrable functions, in H\"{o}lder classes,  in spaces of
generalized functions, in the complex plane, for functions of
several variables and in various other contexts. We could not
find, however, a good reference for fractional
integrodifferentiation of measures. One reason could be is that
the definition of the fractional integral of a measure is fairly
straightforward.  The inversion problem, nonetheless, i.e. the
definition of a fractional derivative, might not be so trivial.
Here we prove two inversion theorems required in the study of
generalized Stieltjes transforms.  They might be useful in some
other contexts as well, for instance in connection with completely
monotonic functions of positive order, see \cite{KP}.

\bigskip

Let us remind the reader that $\M_{\alpha}$ is the positive cone
comprising non-negative Borel measures supported on $[0,\infty]$
and satisfying (\ref{eq:mu-condition}). Clearly,
$\M_{\alpha}\subset\M_{\beta}$ if $\alpha<\beta$.  There is a
natural involution $N_{\alpha}$ defined on $\M_{\alpha}$ by
(\ref{eq:Nalpha-defined}).  For a measure $\mu\in\M_{\alpha}$ we
denote by $F_{\mu}$ its left-continuous distribution function
normalized by $F_{\mu}(0)=0$. The distribution function $F_{\mu}$
defines the measure $\mu$ uniquely except for a possible atom at
infinity which must be specified separately.

\medskip

\noindent\textbf{Definition~A1.}  Let $\mu\in\M_{\alpha}$. The
measure $\nu:=I^{+}_{\eta}\mu$ is called the left-sided
Riemann-Liouville fractional integral of $\mu$ of order $\eta>0$
if
\begin{equation}\label{eq:RLmu0}
\nu(B):=\frac{1}{\Gamma(\eta)}\int\limits_{B\setminus\{\infty\}}dy
\int\limits_{[0,y)}\frac{\mu(du)}{(y-u)^{1-\eta}}+\mu(B\cap\{\infty\})~~\text{for
each Borel set}~B\subset[0,\infty].
\end{equation}
\rem Formula (\ref{eq:RLmu0}) is certainly a straightforward
generalization of the left-sided Riemann-Liouville fractional
integral as given in \cite[(2.2.1),(2.2.2)]{KST} and
\cite[Chapter~2(5.1),(5.3)]{SKM}. One can check that
$\nu(\{0\})=0$ regardless of $\mu(\{0\})$ by computing the limit
$\lim_{y\to{0}}F_{\nu}(y)=0$.

\medskip
\noindent\textbf{Definition~A2.} Let $\mu\in\M_{\alpha}$. The
measure $\tau:=K^{-}_{\alpha,\eta}\mu$ is called the right-sided
Kober-Erdelyi fractional integral of $\mu$ of order $\eta>0$ if
\begin{equation}\label{eq:RLmuInf}
\tau(B):=\mu(B\cap\{0\})+\frac{1}{\Gamma(\eta)}\!\!\!\int\limits_{B\setminus\{0\}}\!\!\!y^{\alpha-1}dy\left\{
\int\limits_{(y,\infty)}\!\!\frac{\mu(du)}{u^{\eta+\alpha-1}(u-y)^{1-\eta}}+\mu_{\infty}\right\}
\end{equation}
for each Borel set $B\subset[0,\infty]$.

\rem Formula (\ref{eq:RLmuInf}) is certainly a straightforward
generalization of the right-sided Kober-Erdelyi fractional
integral as given in  \cite[(2.6.8)]{KST} and \cite[Chapter~4,
(18.6)]{SKM} with a slight change of notation. Condition
(\ref{eq:mu-condition}) ensures that the right hand side of
(\ref{eq:RLmuInf}) exists.  By definition the measure $\tau$ has
no atom at infinity.

\begin{theorem}\label{th:Leftsided}
If $\mu\in\M_{\alpha}$ then
$\nu=I^{+}_{\eta}\mu\in\M_{\alpha+\eta}$ and given $\nu$ we can
recover $\mu$ from
\begin{equation}\label{eq:mubeta-alpha1}
F_{\mu}(y)=\frac{1}{\Gamma(1+n-\eta)}\left(\frac{d}{dy}\right)^n
\int\limits_{[0,y)}\frac{\nu(du)}{(y-u)^{\eta-n}},
\end{equation}
where $n=[\eta]$ and $\mu(\{\infty\})=\nu(\{\infty\})$.
\end{theorem}
\textbf{Proof.}  To show that $\nu\in\M_{\alpha+\eta}$ compute
\begin{multline*}
\int\limits_{[0,\infty)}\frac{\nu(dt)}{(1+t)^{\alpha+\eta}}
=\frac{1}{\Gamma(\eta)}\int\limits_{[0,\infty)}\frac{dt}{(1+t)^{\alpha+\eta}}
\int\limits_{[0,t)}\frac{\mu(du)}{(t-u)^{1-\eta}}
\\
=\frac{1}{\Gamma(\eta)}\int\limits_{[0,\infty)}\mu(du)
\underbrace{\int\limits_{(u,\infty)}\frac{dt}{(t-u)^{1-\eta}(1+t)^{\alpha+\eta}}}_{=(1+u)^{-\alpha}B(\eta,\alpha)}
=\frac{\Gamma(\alpha)}{\Gamma(\alpha+\eta)}\int\limits_{[0,\infty)}\frac{\mu(du)}{(1+u)^{\alpha}}<\infty
\end{multline*}
according to (\ref{eq:mu-condition}).  Here
$$
B(\eta,\alpha)=\frac{\Gamma(\alpha)\Gamma(\eta)}{\Gamma(\alpha+\eta)}
$$
is Euler's beta function.  The interchange of the order of
integrations is legitimate by Tonelli's theorem
\cite[Theorem~13.8]{Schilling}, \cite[Chapter~20, Corollary~7]{RF}
or \cite[Theorem~3.4.5]{Bogachev}.  The proof of
(\ref{eq:mubeta-alpha1}) is a paraphrase of the standard proof of
the inversion formula for the Riemann-Liouville fractional
integral \cite[Theorem~2.4]{SKM} except that we recover the
distribution function of the measure $\mu$ so that we
differentiate one time less the standard and we employ Tonelli's
theorem to justify the interchange of integrations. Substitute
(\ref{eq:RLmu0}) into the integral on the right hand side of
(\ref{eq:mubeta-alpha1}) (recall that $n=[\eta]$):
\begin{multline*}
\int\limits_{[0,y)}\frac{\nu(du)}{(y-u)^{\eta-n}}
\\
=\frac{1}{\Gamma(\eta)}\int\limits_{[0,y)}\frac{du}{(y-u)^{\eta-n}}
\int\limits_{[0,u)}\frac{\mu(dt)}{(u-t)^{1-\eta}}=\frac{1}{\Gamma(\eta)}\int\limits_{[0,y)}\mu(dt)
\int\limits_{(t,yh)}\frac{du}{(u-t)^{1-\eta}(y-u)^{\eta-n}}
\\
=\frac{\Gamma(1-\eta+n)}{n!}\int\limits_{[0,y)}(y-t)^n\mu(dt)=\Gamma(1-\eta+n)\int\limits_{[0,y)}\!\!dt_1
\int\limits_{[0,t_1)}\!\!dt_2\cdots\int\limits_{[0,t_{n-1})}\!\!F_{\mu}(t_{n})dt_{n}.
\end{multline*}
Since $F_{\mu}(t_{n})$ is non-decreasing it is locally Lebesgue
integrable which implies that the function on the right belongs to
$AC^n[0,R]$ (the function and $n-1$ its derivative are absolutely
continuous) for any $R>0$. Hence, we can recover $F_{\mu}$ by
$n$-fold differentiation.~~$\square$

\begin{theorem}\label{th:Rightsided}
If $\mu\in\M_{\alpha}$ then
$\tau:=K^{-}_{\alpha,\eta}(\mu)\in\M_{\alpha+\eta}$ and  given
$\tau$ we can recover $\mu$ from
\begin{multline}\label{eq:rhobeta-alph1}
F_{\mu}(y)=\tau(\{0\})+\frac{1}{\Gamma(n+1-\eta)}\times
\\
\biggl\{\!\alpha\!\!\int\limits_{(0,y)}\!\!x^{\alpha-1}dx
\left(-x^2\frac{d}{dx}\right)^nx^{\eta-n}\!\!\!\!\int\limits_{(x,\infty)}\!\!\!\!\frac{\tau(ds)}{s^{\alpha+n}(s-x)^{\eta-n}}
-y^{\alpha}\left(-y^2\frac{d}{dy}\right)^ny^{\eta-n}\!\!\!\!\int\limits_{(y,\infty)}\!\!\!\!\frac{\tau(ds)}{s^{\alpha+n}(s-y)^{\eta-n}}\biggr\},
\end{multline}
where $n=[\eta]$ and
$\mu_{\infty}=\alpha\Gamma(\eta)\lim\limits_{y\to\infty}y^{-\alpha}F_{\tau}(y)$.
\end{theorem}
\textbf{Proof.} To show that $\tau\in\M_{\alpha+\eta}$ compute
\begin{multline*}
\int\limits_{(0,\infty)}\frac{\tau(dt)}{(1+t)^{\alpha+\eta}}=\frac{1}{\Gamma(\eta)}\int\limits_{[0,\infty)}\frac{t^{\alpha-1}dt}{(1+t)^{\alpha+\eta}}
\int\limits_{(t,\infty)}\frac{\mu(du)}{u^{\eta+\alpha-1}(u-t)^{1-\eta}}
\\
=\frac{1}{\Gamma(\eta)}\int\limits_{(0,\infty)}u^{1-\eta-\alpha}\mu(du)
\underbrace{\int\limits_{(0,u)}\frac{t^{\alpha-1}dt}{(1+t)^{\alpha+\eta}(u-t)^{1-\eta}}}_{=u^{\eta+\alpha-1}(1+u)^{-\alpha}B(\alpha,\eta)}
=\frac{\Gamma(\alpha)}{\Gamma(\alpha+\eta)}\int\limits_{(0,\infty)}\frac{\mu(du)}{(1+u)^{\alpha}}<\infty
\end{multline*}
The interchange of integrations is again justified by Tonelli's
theorem.  To prove formula (\ref{eq:rhobeta-alph1}) assume for the
moment that $\mu_{\infty}=0$ and substitute (\ref{eq:RLmuInf})
into (\ref{eq:rhobeta-alph1}):
\begin{multline}\label{eq:tau-chain}
\alpha\!\!\int\limits_{(0,y)}\!\!x^{\alpha-1}dx
\left(-x^2\frac{d}{dx}\right)^nx^{\eta-n}\!\!\!\!\int\limits_{(x,\infty)}\!\!\!\!\frac{\tau(ds)}{s^{\alpha+n}(s-x)^{\eta-n}}
-y^{\alpha}\left(-y^2\frac{d}{dy}\right)^ny^{\eta-n}\!\!\!\!\int\limits_{(y,\infty)}\!\!\!\!\frac{\tau(ds)}{s^{\alpha+n}(s-y)^{\eta-n}}
\\=\frac{\alpha}{\Gamma(\eta)}\!\!\int\limits_{(0,y)}\!\!x^{\alpha-1}dx
\left(-x^2\frac{d}{dx}\right)^nx^{\eta-n}\!\!\!\!\int\limits_{(x,\infty)}\!\!\!\!\frac{s^{\alpha-1}ds}{s^{\alpha+n}(s-x)^{\eta-n}}
\int\limits_{(s,\infty)}\frac{\mu(du)}{u^{\eta+\alpha-1}(u-s)^{1-\eta}}
\\
-\frac{y^{\alpha}}{\Gamma(\eta)}\left(-y^2\frac{d}{dy}\right)^ny^{\eta-n}\!\!\!\!\int\limits_{(y,\infty)}\!\!\!\!\frac{s^{\alpha-1}ds}{s^{\alpha+n}(s-y)^{\eta-n}}
\int\limits_{(s,\infty)}\frac{\mu(du)}{u^{\eta+\alpha-1}(u-s)^{1-\eta}}
\\
=\frac{\alpha}{\Gamma(\eta)}\!\!\int\limits_{(0,y)}\!\!x^{\alpha-1}dx
\biggl\{\left(\frac{d}{dt}\right)^n(1/t)^{\eta-n}\!\!\!\!\int\limits_{(1/t,\infty)}\!\!\!\!\frac{ds}{s^{n+1}(s-1/t)^{\eta-n}}
\int\limits_{(s,\infty)}\frac{\mu(du)}{u^{\eta+\alpha-1}(u-s)^{1-\eta}}\biggr\}_{\!\!|\,t=1/x}
\\
-\frac{y^{\alpha}}{\Gamma(\eta)}\biggl\{\left(\frac{d}{dt}\right)^n(1/t)^{\eta-n}\!\!\!\!\int\limits_{(1/t,\infty)}\!\!\!\!\frac{ds}{s^{n+1}(s-1/t)^{\eta-n}}
\int\limits_{(s,\infty)}\frac{\mu(du)}{u^{\eta+\alpha-1}(u-s)^{1-\eta}}\biggr\}_{\!\!|\,t=1/y},
\end{multline}
where we have used the formula
$$
\left(-y^2\frac{d}{dy}\right)^n\varphi(y)=\biggl\{\left(\frac{d}{dt}\right)^n\varphi(1/t)\biggr\}_{\!\!|\,t=1/y}.
$$
Further, exchange of the order of integrations justified by
Tonelli's theorem yields:
\begin{multline*}
\int\limits_{(1/t,\infty)}\!\!\!\!\frac{ds}{s^{n+1}(s-1/t)^{\eta-n}}\int\limits_{(s,\infty)}\frac{\mu(du)}{u^{\eta+\alpha-1}(u-s)^{1-\eta}}
\\
=\int\limits_{(1/t,\infty)} \frac{\mu(du)}{u^{\eta+\alpha-1}}
\underbrace{\int\limits_{(1/t,u)}\frac{ds}{s^{n+1}(s-1/t)^{\eta-n}(u-s)^{1-\eta}}}_{=t^{\eta-n}
u^{\eta-n-1}(ut-1)^n B(\eta,1-\eta+n)}
=B(\eta,1-\eta+n)t^{\eta-n}\int\limits_{(1/t,\infty)}\!\!\frac{\mu(du)}{u^{n+\alpha}}(ut-1)^n.
\end{multline*}

We will show now that
$$
\left(\frac{d}{dt}\right)^n\int\limits_{(1/t,\infty)}\!\!\frac{\mu(du)}{u^{n+\alpha}}(ut-1)^n
=n!\int\limits_{(1/t,\infty)}u^{-\alpha}\mu(du).
$$
Denote by $\mu^{*}$ the image of the measure $\mu$ under the
mapping $\lambda(u)=1/u$, so that for each Borel set
$A\subset(0,\infty)$ we have  $\mu^{*}(A):=\mu(\lambda^{-1}(A))$.
Then
\begin{multline*}
\left(\frac{d}{dt}\right)^n\int\limits_{(1/t,\infty)}\!\!\frac{\mu(du)}{u^{n+\alpha}}(ut-1)^n
=\left(\frac{d}{dt}\right)^n\int\limits_{(0,t)}\!\!s^{n+\alpha}(t/s-1)^n\mu^{*}(ds)
=\left(\frac{d}{dt}\right)^n\int\limits_{(0,t)}\!\!(t-s)^ns^{\alpha}\mu^{*}(ds)
\\
=n!\left(\frac{d}{dt}\right)^n\int\limits_{(0,t)}\!\!dt_1\int\limits_{(0,t_1)}\!\!dt_2\cdots
\int\limits_{(0,t_{n})}s^{\alpha}\mu^{*}(ds)=n!\int\limits_{(0,t)}s^{\alpha}\mu^{*}(ds)
=n!\int\limits_{(1/t,\infty)}u^{-\alpha}\mu(du)
\end{multline*}
Finally, by integration by parts for Lebesgue-Stieltjes integral
(see \cite[Theorem 6.2.2]{Carter-vanBrunt}) we obtain:
\begin{multline*}
\alpha\!\!\int\limits_{(0,y)}\!\!x^{\alpha-1}dx
\left(-x^2\frac{d}{dx}\right)^nx^{\eta-n}\!\!\!\!\int\limits_{(x,\infty)}\!\!\!\!\frac{\tau(ds)}{s^{\alpha+n}(s-x)^{\eta-n}}
-y^{\alpha}\left(-y^2\frac{d}{dy}\right)^ny^{\eta-n}\!\!\!\!\int\limits_{(y,\infty)}\!\!\!\!\frac{\tau(ds)}{s^{\alpha+n}(s-y)^{\eta-n}}
\\
=\Gamma(1-\eta+n)\biggl(\alpha\!\!\int\limits_{(0,y)}\!\!x^{\alpha-1}dx\int\limits_{(x,\infty)}u^{-\alpha}\mu(du)
-y^{\alpha}\int\limits_{(y,\infty)}u^{-\alpha}\mu(du)\biggr)
\\
=\Gamma(1-\eta+n)\biggl(x^{\alpha}\int\limits_{(x,\infty)}\!\!\!u^{-\alpha}\mu(du)\biggl|_{0}^{y}
+\int\limits_{(0,y)}\mu(dx)
-y^{\alpha}\!\!\!\int\limits_{(y,\infty)}u^{-\alpha}\mu(du)\biggr)
\\
=\Gamma(1-\eta+n)\int\limits_{(0,y)}\mu(dx)
=\Gamma(1-\eta+n)(F_{\mu}(y)-F_{\mu}(0+))=\Gamma(1-\eta+n)(F_{\mu}(y)-\mu(\{0\})).
\end{multline*}
Since $\mu(\{0\})=\tau(\{0\})$ by (\ref{eq:RLmuInf}) this proves
formula (\ref{eq:rhobeta-alph1}). While integrating by parts we
have also used the following limit:
\begin{multline*}
\lim\limits_{x\to{0}}x^{\alpha}\int\limits_{(x,\infty)}\!\!\!u^{-\alpha}\mu(du)
=\lim\limits_{x\to{0}}x^{\alpha}\int\limits_{(x,1)}\!\!\!u^{-\alpha}dF_{\mu}(u)
=\lim\limits_{x\to{0}}x^{\alpha}\biggl(u^{-\alpha}F_{\mu}(u)\biggl|^{1}_{x}+\alpha
\int\limits_{(x,1)}\!\!\!u^{-\alpha-1}F_{\mu}(u)du\biggr)
\\
=\lim\limits_{x\to{0}}\biggl({\alpha}x^{\alpha}\int\limits_{(x,1)}\!\!\!u^{-\alpha-1}F_{\mu}(u)du
-F_{\mu}(x)\biggr)
=\lim\limits_{x\to{0}}\frac{-x^{-\alpha-1}F_{\mu}(x)}{(1/\alpha)(-\alpha)x^{-\alpha-1}}
-F_{\mu}(0+)=0.
\end{multline*}
Here the first equality is due to (\ref{eq:mu-condition}), the
second is integration by parts and the preultimate is
L'H\^{o}pital's rule applied if
$\int_{(x,1)}\!u^{-\alpha-1}F_{\mu}(u)du$ is unbounded.

If $\mu_{\infty}\ne{0}$ we need to add the following term to the
third line of (\ref{eq:tau-chain}):
\begin{multline*}
\frac{\alpha\mu_{\infty}}{\Gamma(\eta)}\!\!\int\limits_{(0,y)}\!\!x^{\alpha-1}dx
\left(\!\!-x^2\frac{d}{dx}\right)^nx^{\eta-n}\!\!\!\!\underbrace{\int\limits_{(x,\infty)}\!\!\!\!\frac{s^{\alpha-1}ds}{s^{\alpha+n}(s-x)^{\eta-n}}}_{=B(\eta,1-\eta+n)x^{-\eta}}
-\frac{\mu_{\infty}y^{\alpha}}{\Gamma(\eta)}\left(\!\!-y^2\frac{d}{dy}\right)^ny^{\eta-n}\!\!\!\!
\overbrace{\int\limits_{(y,\infty)}\!\!\!\!\frac{s^{\alpha-1}ds}{s^{\alpha+n}(s-y)^{\eta-n}}}^{{=B(\eta,1-\eta+n)y^{-\eta}}}
\\
=\frac{\alpha\mu_{\infty}}{\Gamma(\eta)}\!\!\int\limits_{(0,y)}\!\!x^{\alpha-1}dx
\left(\!\!-x^2\frac{d}{dx}\right)^n\!\!x^{-n}-\frac{\mu_{\infty}y^{\alpha}}{\Gamma(\eta)}\left(\!\!-y^2\frac{d}{dy}\right)^n\!\!y^{-n}
=\frac{\alpha\mu_{\infty}n!}{\Gamma(\eta)}\!\!\int\limits_{(0,y)}\!\!x^{\alpha-1}dx
-\frac{\mu_{\infty}n!y^{\alpha}}{\Gamma(\eta)}=0,
\end{multline*}
where
$$
\left(\!\!-x^2\frac{d}{dx}\right)^n\!\!x^{-n}=\biggl\{\left(\frac{d}{dt}\right)^nt^n\biggr\}_{\!\!|\,t=1/x}=n!.
$$
This shows that (\ref{eq:tau-chain}) is still valid for
$\mu_{\infty}\ne{0}$.

Finally, to recover the atom at infinity $\mu_{\infty}$ we compute
\begin{multline*}
\lim\limits_{x\to\infty}x^{-\alpha}F_{\tau}(x)
=\lim\limits_{x\to\infty}[\mu(\{0\})x^{-\alpha}]+\frac{1}{\Gamma(\eta)}
\lim\limits_{x\to\infty}x^{-\alpha}\!\int\limits_{(0,x)}\!\!\!y^{\alpha-1}dy\int\limits_{(y,\infty)}\!\!\frac{\mu(du)}{u^{\eta+\alpha-1}(u-y)^{1-\eta}}
\\
 +\frac{\mu_{\infty}}{\Gamma(\eta)}
\lim\limits_{x\to\infty}x^{-\alpha}\!\int\limits_{(0,x)}\!\!\!y^{\alpha-1}dy=\frac{\mu_{\infty}}{\alpha\Gamma(\eta)}
+\frac{1}{\Gamma(\eta)}
\lim\limits_{x\to\infty}x^{-\alpha}\!\int\limits_{(0,x)}\!\!\!y^{\alpha-1}dy\int\limits_{(y,\infty)}\!\!\frac{\mu(du)}{u^{\eta+\alpha-1}(u-y)^{1-\eta}}.
\end{multline*}
In order to show that the last limit is zero we interchange the
order of integrations (justified again by Tonelli's theorem):
\begin{multline*}
\int\limits_{(0,x)}\!\!\!y^{\alpha-1}dy\int\limits_{(y,\infty)}\!\!\frac{\mu(du)}{u^{\eta+\alpha-1}(u-y)^{1-\eta}}
=\int\limits_{(0,x)}u^{1-\eta-\alpha}\mu(du)\overbrace{\int\limits_{(0,u)}(u-y)^{\eta-1}y^{\alpha-1}dy}^{=B(\alpha,\eta)u^{\alpha+\eta-1}}
\\
+\int\limits_{[x,\infty)}u^{1-\eta-\alpha}\mu(du)\!\!\!\!\underbrace{\int\limits_{(0,x)}(u-y)^{\eta-1}y^{\alpha-1}dy}_{=(1/\alpha)x^{\alpha}u^{\eta-1}{_2F_1}(\alpha,1-\eta;1+\alpha;x/u)}
\!\!\!\!\!=F_{\mu}(x)-\mu(\{0\})
\\
+\frac{1}{\alpha}x^{\alpha}\int\limits_{[x,\infty)}u^{-\alpha}{_2F_1}(\alpha,1-\eta;1+\alpha;x/u)\mu(du),
\end{multline*}
where ${_2F_1}$ is the Gauss hypergeometric function \cite[Chapter
2]{AAR}.  Hence we need to prove that:
$$
\lim\limits_{x\to\infty}x^{-\alpha}F_{\mu}(x)=0
$$
and
$$
\lim\limits_{x\to\infty}\int\limits_{[x,\infty)}u^{-\alpha}{_2F_1}(\alpha,1-\eta;1+\alpha;x/u)\mu(du)=0.
$$
Both equalities follow from (\ref{eq:mu-condition}): the first was
proved by Widder \cite[Corollary 3a.3]{Widder}, the second follows
from the fact that ${_2F_1}(a,b;c;x)$ is bounded on $[0,1]$ if
$c>a+b$ by the Gauss formula \cite[Theorem 2.2.2]{AAR}.
~~$\square$

\rem The operator $(-x^2D_x)^n$ encountered in
(\ref{eq:rhobeta-alph1}) can be expanded as follows
$$
(-x^2D_x)^nf=\sum\limits_{m=1}^{n}a(n,m)x^{n+m}f^{(m)}(x),
$$
where the numbers
$$
a(n,m)=(-1)^n\frac{n!}{m!}\binom{n-1}{m-1}
$$
are known as Lah numbers \cite[A008297]{OEIS} satisfying
$$
a(n+1,m)=(n+m)a(n,m)+a(n,m-1).
$$
Applying $-x^2D_x$ to the above expansion we see the same
recurrence which given the same initial values furnishes a proof
of the above expansion.

\rem Another way to obtain a representation for $\mu$ via $\tau$
for $\eta>1$ is the following.  Denote
$\tilde{\mu}=K^{-}_{\alpha,1}(\mu)$, i.e. according to
(\ref{eq:RLmuInf})
$$
\tilde{\mu}(dx)=\mu(\{0\})\delta_0+x^{\alpha-1}dx\int\limits_{(x,\infty)}u^{-\alpha}\mu(du)=\mu(\{0\})\delta_0+\phi_{\tilde{\mu}}(x)dx,
$$
where the last equality is the definition of
$\phi_{\tilde{\mu}}(x)$. It is easy to verify that this formula is
inverted as follows
$$
F_{\mu}(y-0)-F_{\mu}(0+)=\alpha\int\limits_{(0,y)}\tilde{\mu}(dx)-y\phi_{\tilde{\mu}}(y).
$$
  If $\eta>1$ we have according to
(\ref{eq:RLmuInf}) for $\tilde{\mu}\in\M_{\alpha+1}$ and the
semigroup property of Kober-Erdeliy operator \cite[(2.6.24)]{KST}
\begin{multline}\label{eq:phibeta}
\tau(dx):=K^{-}_{\alpha,\eta}(\mu)=K^{-}_{\alpha+1,\eta-1}(\tilde{\mu})
\\
=\mu(\{0\})\delta_0+\frac{x^{\alpha}dx}{\Gamma(\eta-1)}
\int\limits_{(x,\infty)}\frac{\phi_{\tilde{\mu}}(u)du}{u^{\alpha+\eta-1}(u-x)^{2-\eta}}=\mu(\{0\})\delta_0+\phi_{\tau}(x)dx.
\end{multline}
Here we have the standard Riemann-Liouvile fractional integral of
the function $u^{1-\alpha-\eta}\phi_{\tilde{\mu}}(u)$. We cannot,
however, use the Riemann-Liouvile  fractional derivative to invert
the above formula, since, in general the integral in its
definition will diverge. Instead, we can employ Caputo's
fractional derivative \cite[section~2.4]{KST} to invert
(\ref{eq:phibeta}):
$$
u^{1-\alpha-\eta}\phi_{\tilde{\mu}}(u)={}^C\!\!D_{-}^{\eta-1}[x^{-\alpha}\phi_{\tau}(x)](u)
=\frac{(-1)^n}{\Gamma(n-\eta+1)}
\int\limits_{(u,\infty)}\frac{[x^{-\alpha}\phi_{\tau}(x)]^{(n)}dx}{(x-u)^{\eta-n}},
$$
where $n=[\eta]$ and
$$
\tilde{\mu}(du)=\mu(\{0\})\delta_0+\frac{(-1)^nu^{\alpha+\eta-1}du}{\Gamma(n-\eta+1)}
\int\limits_{(u,\infty)}\frac{[x^{-\alpha}\phi_{\tau}(x)]^{(n)}dx}{(x-u)^{\eta-n}},
$$
so that
\begin{multline*}
F_{\mu}(y)=\mu(\{0\})\delta_0+\frac{(-1)^n\alpha}{\Gamma(n-\eta+1)}\int\limits_{(0,y)}u^{\alpha+\eta-1}du
\int\limits_{(u,\infty)}\frac{[x^{-\alpha}\phi_{\tau}(x)]^{(n)}}{(x-u)^{\eta-n}}dx
\\
-\frac{(-1)^ny^{\alpha+\eta}}{\Gamma(n-\eta+1)}
\int\limits_{(y,\infty)}\frac{[x^{-\alpha}\phi_{\tau}(x)]^{(n)}}{(x-y)^{\eta-n}}dx.
\end{multline*}
This formula is, however, less general then
(\ref{eq:rhobeta-alph1}) since it requires that
$\phi_{\tau}(x)\in{AC^n(0,\infty)}$ which is not guaranteed by
$\mu\in\M_{\alpha}$.  For $0<\eta<1$ both formulas take the same
form.

The relation between $I^{+}_{\eta}$ and $K^{-}_{\alpha,\eta}$ is
revealed in the following theorem.
\begin{theorem}
Suppose $\mu\in\M_{\alpha}$. Then
$N_{\alpha+\eta}I^{+}_{\eta}\mu=K^{-}_{\alpha,\eta}N_{\alpha}\mu$.
\end{theorem}
\textbf{Proof.} Indeed, if $f$ is given by (\ref{eq:mu-rep}) and
$\nu=I^{+}_{\eta}\mu$ then by Theorem~\ref{lm:3}
$$
f(z)-\mu_{\infty}=\int\limits_{[0,\infty)}\frac{\mu(du)}{(u+z)^{\alpha}}=\frac{\Gamma(\alpha+\eta)}{\Gamma(\alpha)}
\int\limits_{(0,\infty)}\frac{\nu(du)}{(u+z)^{\alpha+\eta}}=\frac{\Gamma(\alpha+\eta)}{\Gamma(\alpha)}\int\limits_{(0,\infty)}\frac{\tau_1(dt)}{(1+tz)^{\alpha+\eta}},
$$
where $\tau_1=N_{\alpha+\eta}\nu=N_{\alpha+\eta}I^{+}_{\eta}\mu$.
On the other hand, if $\rho=N_{\alpha}\mu$ then according to the
comment after formula (\ref{eq:Nalpha-defined}) and by
Theorem~\ref{lm:31} we have
$$
f(z)-\mu_{\infty}=\int\limits_{(0,\infty)}\frac{\rho(dt)}{(1+tz)^{\alpha}}+\frac{\mu_0}{z^{\alpha}}
=\frac{\Gamma(\alpha+\eta)}{\Gamma(\alpha)}\int\limits_{(0,\infty)}\frac{\tau_2(dt)}{(1+tz)^{\alpha+\eta}},
$$
where
$\tau_2=K^{-}_{\alpha,\eta}\rho=K^{-}_{\alpha,\eta}N_{\alpha}\mu$.
Comparing these two formulas we conclude that $\tau_1=\tau_2$ due
to uniqueness of the representing measure.~~$\square$

\paragraph{5. Acknowledgements.} We thank Christian Berg, Leonid Kovalev,
Jos\'{e} Luis L\'{o}pez and Sergei Sitnik for numerous useful
discussions and Alex Gomilko for counterexample in Remark~7. We
acknowledge the financial support of the Russian Basic Research
Fund (grant 11-01-00038-a) and the Far Eastern Branch of the
Russian Academy of Sciences (grant 12-III-A-01-007).

\end{document}